\documentclass{itogi20171}

\currentyear{}
\currentvolume{}

\usepackage[russian]{babel}
\usepackage[cp1251]{inputenc}

\usepackage{graphicx}

\usepackage{enumerate} 

\usepackage{wasysym} 

\usepackage{xcolor} 

\usepackage[
unicode,colorlinks,linkcolor=blue,citecolor=red,bookmarksopen,pdfhighlight=/N]{hyperref}

\newtheorem{thm}{Теорема}[section]

\newtheorem*{thMR}{Критерий Мальявена\,--\,Рубела}
\newtheorem{proposition}{Предложение}[section]
\newtheorem{corollary}{Следствие}[section]

\theoremstyle{definition}

\newtheorem{remark}{Замечание}[section]
\newtheorem{example}{Пример}[section]

\renewcommand{\leq}{\leqslant} 
\renewcommand{\geq}{\geqslant}
\newcommand{\RR}{\mathbb{R}} 
\newcommand{\CC}{\mathbb{C}} 
\newcommand{\NN}{\mathbb{N}} 
\newcommand{\ZZ}{\mathbb{Z}}

\DeclareMathOperator{\Zero}{Zero} 
 
\DeclareMathOperator{\strip}{str} 
\DeclareMathOperator{\dd}{d}
\DeclareMathOperator{\intr}{\mathsf{int}}
\DeclareMathOperator{\Hol}{\mathsf{Hol}}
\DeclareMathOperator{\dens}{\mathsf{dens}}
\DeclareMathOperator{\Exp}{\mathsf{Exp}}
\DeclareMathOperator{\diam}{\mathsf{diam}}
\DeclareMathOperator{\width}{\mathsf{breadth}}
\DeclareMathOperator{\spf}{\mathsf{spf}}
\DeclareMathOperator{\Spf}{\mathsf{Spf}}
\DeclareMathOperator{\breadth}{\mathsf{width}}

\renewcommand{\Re}{\operatorname{Re}}
\renewcommand{\Im}{\operatorname{Im}}

\begin{document}

\title[Полнота экспоненциальной системы в геометрических терминах]{Полнота экспоненциальной системы 
в геометрических терминах  ширины, широты и диаметра}

\author[Хабибуллин Булат Нурмиевич]{Б. Н.  Хабибуллин}
\address{Башкирский государственный университет, 
Институт математики с ВЦ УФИЦ РАН}
\email{khabib-bulat@mail.ru}
\thanks{Работа выполнена в рамках государственного задания Министерства науки и высшего образования 
Российской Федерации (код научной темы FMRS-2022-0124) при поддержке   Министерства просвещения Российской Федерации в рамках государственного задания (соглашение № 073-03-2023-010 от 26.01.2023).}

\author[Кудашева Елена Геннадьевна]{Е. Г. Кудашева}
\address{Башкирский государственный педагогический университет им. М.~Акмуллы}
\email{lena\_kudasheva@mail.ru}

\author[Салимова Анна  Евгениевна]{А. Е. Салимова}
\address{Башкирский государственный университет}
\email{anegorova94@bk.ru}




\UDC{517.538.2, 517.547.22,  514.17, 517.574}

\begin{abstract}
Устанавливается критерий полноты экспоненциальной системы в пространствах функций, непрерывных  на выпуклом компакте и голоморфных во внутренности этого компакта, а также в пространствах голоморфных функций в выпуклой области в терминах ширины компакта или области в направлении. Основные результаты формулируются исключительно через соотношения между шириной в направлении, широтой или диметром компакта или области с одной стороны и так называемыми  
логарифмическими субмерами или  логарифмическими блок-плотностями  распределения показателей экспоненциальной системы с другой.

\end{abstract}


\maketitle

\tableofcontents

\section{Введение. Формулировки критериев в терминах ширины множества}

Одноточечные множества $\{a\}$ часто записываем без фигурных скобок, т.е. просто как $a$. Так,   $\NN_0:=0\bigcup \NN=\{0,1, 2, \dots\}$ для  множества $\mathbb N:=\{1,2, \dots\}$ {\it натуральных чисел,\/}  а множество 
$\overline \NN_0:=\NN_0\bigcup +\infty$ --- верхнее порядковое пополнение  множества $\NN_0$ со стандартным отношением порядка $\leq$ точной верхней гранью $+\infty:=\sup \NN_0\notin \NN_0$, для которой выполнены неравенства  $n\leq  +\infty$ при всех $n\in \overline \NN_0$. \textit{Множество всех действительных чисел\/} $\RR$ с таким же отношением порядка $\leq$  рассматриваем и  как \textit{вещественную ось\/} \textit{в комплексной плоскости\/} $\CC$,  а \textit{множество всех положительных чисел\/} $\RR^+:=\bigl\{x\in \RR\bigm| 0\leq x\bigr\}$ как 
\textit{положительную полуось\/} в $\RR$ или в $\CC$. Порядковое пополнение множества $\RR$ верхней гранью $+\infty  :=\sup \RR$ и нижней  гранью $-\infty  :=\inf \RR$ даёт \textit{расширенные\/} вещественную ось  $\overline \RR:=+\infty \bigcup \RR\bigcup -\infty$ и положительную полуось $\overline \RR^+:=\RR^+\bigcup +\infty$. 
Для \textit{пустого множества\/} $\emptyset$  по определению $\sup \emptyset:=-\infty$ и $\inf \emptyset :=+\infty$.

Система векторов из  топологического векторного пространства \textit{полна\/} в нём,  
если замыкание линейной оболочки этой системы совпадает с этим пространством. Здесь мы обсуждаем только полноту экспоненциальных систем 
в  функциональных пространствах. Истоки нашего исследования в следующем результате П. Мальявена и Л.\,А. Рубела из их совместной статьи \cite{MR} начала 1960-х гг. Для формулировки их критерия  \textit{открытую\/} и \textit{ замкнутую полосы ширины $b\in \overline \RR^+$ в $\CC$, 
 симметричные\/} относительно вещественной оси\/ $\RR$ и нуля, обозначаем соответственно через 
\begin{equation}\label{strip}
\strip_{b/2}:=\bigl\{ z\in \CC \colon |\Im z|<b/2\bigr\},\quad
\overline{\strip}_{b/2}:=\bigl\{ z\in \CC \colon |\Im z|\leq b/2\bigr\}.
\end{equation}

\begin{thMR}[{\cite[теорема 9.1]{MR}}] Пусть $0<b\in \RR^+$ и $Z=(z_n)_{n\in \NN}$ --- возрастающая последовательность  положительных попарно различных чисел, для которой  последовательность 
$(n/z_n)_{n\in \NN}$ ограничена, т.е. $Z$ конечной верхней плотности. 
Экспоненциальная система  $\Exp^{Z}$, состоящая из функций  $z\underset{z\in \CC}{\longmapsto} e^{z_nz}$, $n\in \NN$,  полна в пространстве   непрерывных на $\overline{\strip}_{b/2}$ и голоморфных на  ${\strip}_{b/2}$ 
функций, снабжённом  топологией равномерной сходимости на компактах из $\overline{\strip}_{b/2}$, если и только если \underline{не} существует числа $C\in \RR$, для которого 
\begin{equation}\label{sMR}
\sum_{r<z_n\leq R}\frac{1}{z_n}\leq \frac{b}{2\pi}\ln\frac{R}{r}+C \quad\text{при всех $0<r<R<+\infty$.}
\end{equation} 
\end{thMR}
\begin{remark}\label{rem1}
Оригинальная формулировка \cite[теорема 9.1]{MR} --- это  \textit{критерий \underline{не}полноты экспоненциальной системы\/} $\Exp^{-Z}$ из функций  $z\underset{z\in \CC}{\longmapsto} e^{-z_nz}$, $n\in \NN$
для полос \eqref{strip} ширины  $2\pi b$ вместо $b$. Она   эквивалентна  сформулированному нами критерию полноты Мальявена\,--\,Рубела.  
\end{remark}
Для компакта $K$ в $\CC$  через $C(K)$ обозначаем \textit{банахово пространство-алгебру непрерывных функций $f\colon K\to \CC$ с\/  $\sup$-нормой} 
\begin{equation}\label{CKn}
\|f\|_{C(K)}:=\sup\Bigl\{\bigl|f(z)\bigr|\Bigm| z\in K \Bigr\}.
\end{equation}
Для открытого подмножества $O\subset \CC$ через $\Hol(O)$ обозначаем пространство голоморфных функций $f\colon O\to \CC$ с топологией равномерной сходимости на всех компактах $K\subset O$, определяемой полунормами 
\eqref{CKn}. Для компакта $K\subset \CC$ с \textit{внутренностью\/} $\intr K$ через $C(K)\bigcap \Hol(\intr K)$  обозначаем банахово пространство непрерывных на $K$ и голоморфных на внутренности $\intr K$ функций
$f\colon K\to \CC$ с $\sup$-нормой  \eqref{CKn}. Очевидно,  если $\intr K=\emptyset$, то $C(K)\bigcap \Hol(\intr K)=C(K)$.

Всюду далее через $Z$ обозначаем \textit{распределение точек\/} на комплексной плоскости\/ $\mathbb C$, среди которых  могут быть повторяющиеся  и, вообще говоря, даже бесконечное  количество раз. Распределение точек $Z$ однозначно определяется функцией, действующей из $\mathbb C$ в $\overline \NN_0$ и равной в каждой точке $z\in \mathbb C$ количеству  повторений этой точки $z$  в распределение точек $Z$. Для такой функции, которую часто называют \textit{функцией кратности,\/} или \textit{дивизором\/} распределения точек  $Z$ \cite[пп.~0.1.2--0.1.3]{Khsur}, 
сохраняем то же обозначение $Z$. Другими словами, $Z(z)$ --- это количество вхождений   точки  $z\in \CC$ в $Z$ и  
пишем $z\in Z$,  если $Z(z)>0$.
Два распределения точек $Z$ и $W$ \textit{совпадают\/} и пишем $Z=W$, если $Z(z)\underset{z\in \CC}{\equiv}W(z)$.  
Пишем $Z\subset W$, если $Z(z)\leq W(z)$ для всех $z\in \CC$. \textit{Объединение \/} $Z\bigcup W$ определяется 
тождеством $(Z\bigcup W)(z)\underset{z\in \CC}{\equiv} Z(z)+ W(z)$, а \textit{разность} $Z\setminus W $
при условии $W\subset Z$ --- тождеством $(Z\setminus W)(z)\underset{z\in \CC}{\equiv} Z(z)- W(z)$. Распределение точек  $Z$ можно эквивалентным образом трактовать  и как меру со значениями в  $\overline \NN_0$ с тем же обозначением  
\begin{equation}\label{Z}
Z(S):=\sum_{z\in S} Z(z)\in \overline \NN_0\quad\text{\it для любого  $S\subset \CC$.}
\end{equation}
При такой трактовке пересечение $Z\bigcap S$ однозначно определяется сужением меры $Z$ на $S$,  
а  для положительной функции $f\colon S\to \overline \RR^+$ можно корректно  определить суммы 
\begin{equation}\label{sumZ}
\sum_{\stackrel{z\in Z}{z\in S}}f(z):=\int_Sf\dd Z=:\sum_{z\in Z\bigcap S} f(z)\in \overline \RR^+.
\end{equation}
Произведение числа $w\in \CC\setminus 0$ на  распределению точек $Z$ определяет   распределение точек $wZ$ через 
функцию кратности  $(wZ)(z)\underset{z\in \CC}{\equiv}Z(z/w)$. В частности, $-Z:=(-1)Z$. 
Для числа $\theta \in \RR$ распределение точек $e^{i\theta}Z$ с функцией кратности $(e^{i\theta}Z)(z)\underset{z\in \CC}{\equiv}Z(e^{-i\theta}z)$ называем  \textit{поворотом распределения точек $Z$ на угол $\theta$}.

При $z\in \CC$ и $r\in \RR^+$ через $D_z(r):=\bigl\{z'\in \CC\bigm| |z'-z|<r\bigr\}$ и $\overline D_z(r):=\bigl\{z'\in \CC\bigm| |z'-z|\leq r\bigr\}$, а также 
$\partial \overline D_z(r):=\overline D_z(r)\setminus  D_z(r)$ 
 обозначаем соответственно \textit{открытый и замкнутый круги,}  а также \textit{окружность с центром $z\in \CC$ радиуса $r$.} 

Простейшая характеристика распределение точек $Z$ на $\CC$ 
--- это его \textit{верхняя $p$-плотность,\/}  определяемая в обозначении \eqref{Z} как величина 
\begin{equation}\label{dns}
p\text{-}\overline \dens (Z):=\limsup_{r\to +\infty}  \frac{ Z\bigl(\overline D_0(r)\bigr)}{r^p}\in \overline \RR^+.
\end{equation}
Далее при значении $p=1$ в \eqref{dns} порядок $p$ не упоминаем и пишем просто $\overline \dens (Z):=1\text{-}\overline \dens (Z)$, 
в случае  $\overline \dens (Z)\in \RR^+$ называем $Z$ распределением точек \textit{конечной верхней плотности,} 
а при $\overline \dens (Z)= +\infty$ распределением точек \textit{бесконечной верхней плотности,}

Левая часть неравенства \eqref{sMR}, использованная в  \cite{MR}, \cite[гл.~22]{RC} в качестве  основной 
так называемой логарифмической характеристики распределения точек $Z$ на $\RR^+$, была   распространена  на произвольные распределения точек $Z$ в $\CC$ первоначально  в статьях Б.\,Н.~Хабибуллина \cite{Kha88MZ} и \cite{Kha88}  1988 г. с  дальнейшим развитием и применениями  в его же  работах \cite{Kha89}, \cite{Kha91} 1989-91 гг.  и в последние годы в совместных статьях О.\,А.~Кривошеевой, А.\,Ф.~Кужаева, А.\,И.~Рафикова \cite{KuzRafKri17}, \cite{KriKuz17} 2017 г., а также   А.\,Е.~Салимовой и Б.\,Н.~Хабибуллина \cite{SalKha20}, \cite{SalKha20MZ}, \cite{SalKha21} 2020-21 гг.
Обзор по состоянию тематики вплоть до 2012 г. содержится в  книге Б.\,Н.~Хабибуллина \cite[п.~3.2]{Khsur}. 
Определим   развития упомянутой  логарифмической характеристики   распределения точек $Z$, расположенного произвольным образом на $\CC$. 

Для  величины $a\in \overline \RR$ и функции $f\colon X\to \overline \RR$ их  \textit{положительные части\/}
обозначаем соответственно через $a^+:=\max\{a,0\}$ и $f^+\colon x\underset{x\in X}{\longmapsto} \bigl(f(x)\bigr)^+$.
К примеру, $\Re^+z=\max\bigl\{\Re z, 0\bigr\}$ --- положительная часть действительной части $\Re z$.
Произвольное число $\theta \in \RR$ будем трактовать и  как \textit{направление,\/} определяемое единичным радиус-вектором числа $e^{i\theta}\in \partial \overline D_0(1)$. При такой трактовке  для направления $\theta$ направления $\theta\pm\pi$ \textit{противоположные,\/} направления $\theta\pm\pi/2$ \textit{ортогональные,\/}  направление $-\theta$ \textit{сопряжённое,\/}  а направление $\pi/2-\theta$ \textit{ортогональное сопряжённому.}

\textit{Правую логарифмическую меру интервала $(r,R]\subset \RR^+$  для распределения точек $Z$ на\/ $\CC$} определяем как   величину 
\begin{equation}\label{ell}
\ell_{Z}(r,R)\overset{\eqref{sumZ}}{:=} \sum_{\stackrel{z\in Z}{r<|z|\leq R}}
\Re^+\frac{1}{z} \overset{\eqref{sumZ}}{=}
\sum_{z\in Z\bigcap\bigl(\overline D_0(R)\setminus \overline D_0(r)\bigr)}
\Re^+\frac{1}{z}\in \overline \RR^+,
\end{equation}
 \textit{левая логарифмическую мера интервала $(r,R]\subset \RR^+$  для $Z$ на\/ $\CC$}
--- это величина $\ell_{-Z}(r,R)$ из \eqref{ell} для $-Z$, 
 а \textit{логарифмическую субмеру  интервала $(r,R]\subset \RR^+\setminus 0$
 для  $Z$ на $\CC$} задаём как 
\begin{equation}\label{leZ}
\mathcal{L}_{Z}(r,R)
:=\max\Bigl\{\ell_{Z}(r,R), \ell_{-Z}(r,R)\Bigr\}\in \overline \RR^+.
\end{equation}

 \textit{Ширина подмножества $S\subset \CC$ в направлениях  $\theta\in \RR$}  \cite[{\bf 33}]{BF}, \cite[4.1.1]{Had}, \cite[гл.~I, \S~4]{Sant}, \cite[п.~3.2 ]{Kha14II}, \cite[п.~3.2]{Khsur} --- это функция 
\begin{equation}\label{width}
\width_S\colon \theta \underset{\theta\in \RR}{\longmapsto} \sup\Bigl\{\Re \bigl((z-w)e^{-i\theta}\bigr)\Bigm| z\in S, \; w\in S\Bigr\}\in  \overline \RR^+\bigcup -\infty, 
\end{equation}
которая, очевидно, $\pi$-периодическая, равна $-\infty$ лишь при пустом $S=\emptyset$, а при $S\neq \emptyset$
 --- величина из $\overline \RR^+$,
равная наименьшей ширине замкнутых полос в $\CC$, содержащих $S$ и какую нибудь прямую,   ортогональную  прямой  $\bigl\{te^{i\theta}\in \CC\bigm| t\in \RR\bigr\}$. 
\begin{example} Ширина любого круга в любом направлении равна его диаметру.  
Полосы \eqref{strip}  ширины $b$ в направлениях $\pi/2+\pi k$ при любом целом $k\in \ZZ:=\NN_0\bigcup (-\NN)$, а 
 в любом другом  направлении ширины $+\infty$.   
Для $S\subset \CC$ и числа $a\in \CC$ полагаем $aS:=\bigl\{az\bigm| z\in S\bigr\}\subset \CC$. При таком обозначении   $e^{i(\theta-\pi/2)}\strip_{b/2}$  --- полоса ширины $b$ в направлении $\theta$. 
\end{example}

\begin{thm}\label{th1} Для любых числа  $b\in \RR^+\setminus 0$, направления $\theta\in \RR$ и произвольного  распределения точек $Z$ на\/ $\CC$ следующие три утверждения равносильны:
\begin{enumerate}[{\rm I.}]

\item\label{KhI+}  
Для любой выпуклой области $D\subset \CC$ ширины $\width_D(\theta)\leq b$ в направлении $\theta$
 система 
$$
\Exp^Z:=\Bigl\{w\underset{w\in \CC}{\longmapsto} w^pe^{zw}\Bigm| z\in Z, \; Z(z)-1 \geq p\in \NN_0 \Bigr\}
$$ 
полна в пространстве $\Hol(D)$,  с топологией равномерной сходимости на компактах.

\item\label{KhI}  Для любого выпуклого компакта $K$ ширины $\width_K(\theta)<b$ в направлении $\theta$  
экспоненциальная система  $\Exp^Z$ полна в банаховом пространстве  $C(K)\bigcap \Hol(\intr K)$ 
с нормой \eqref{CKn}. 

\item\label{KhII} Распределение точек $Z$ бесконечной верхней плотности с $\overline \dens (Z)=+\infty$ или же логарифмическая субмера $\mathcal{L}_{e^{i(\pi/2-\theta)}Z}(r,R)$ для поворота  $e^{i(\pi/2-\theta)}Z$ распределения точек $Z$  на угол $\pi/2-\theta$  удовлетворяет неравенству 
\begin{equation}\label{lndZ}
\limsup_{1<s\to +\infty} \frac{1}{\ln s} \limsup_{r\to+\infty} \mathcal{L}_{e^{i(\pi/2-\theta)}Z}(r,sr)\geq \frac{b}{2\pi}.
\end{equation}
\item\label{KhIII} Система $\Exp^Z$ полна в пространстве  $\Hol\bigl(e^{i(\theta-\pi/2)}\strip_{b/2}\bigr)$.
\end{enumerate} 
\end{thm}
Если  при $N\subset \NN$ и некоторой нумерации  $ Z=(z_n)_{n\in N}$   распределения точек $Z$ на $\CC$, 
в которой  каждое число  $z\in \CC$ встречается ровно $ Z(z)$ раз, 
можно так  подобрать {\it последовательность $(m_n)_{n\in N}$  попарно различных целых чисел\/} $m_n\in \ZZ$ и 
число $c\in \mathbb{R}^+$,  что  
\begin{equation}\label{RdZ}
\sum_{n\in N}\Bigl|\frac{1}{z_n}-\frac{c}{m_n}\Bigr|<+\infty, 
\end{equation}
то   {\it внешняя плотность Редхеффера\/} распределения точек  ${\mathrm Z}$ на $\CC$
конечна и не превышает числа $c$ \cite{Red77}, \cite{Kra89}, \cite{Kha94}, \cite[2.1.1]{Khsur}, а сама она равна точной нижней грани  таких  $c\in \RR^+$. 

\begin{example}\label{ex2}
Если распределение \textit{попарно различных точек\/} $Z$ \textit{разделённое\/} в том смысле, что $\inf \bigl\{|z-z'|\bigm| z\in Z, z'\in Z, z\neq z'\bigr\}>0$ и  целиком лежит  в какой-нибудь полосе \eqref{strip}  конечной	 ширины $b\in \RR^+$,
то $Z$ конечной  внешней плотности Редхеффера.
\end{example}
Распределение точек $Z$ на $\CC$ \textit{конечной внешней плотности Редхеффера  вблизи направления\/ $\theta$,}  если найдётся  число $a\in (0,\pi/2]$, для которого  часть 
\begin{equation}\label{Zat}
Z_a^{\theta}:=\Bigl\{z\in Z\Bigm| |\arg z-\theta |< a\Bigr\}
\end{equation}
распределения точек $Z$  после её поворота $e^{-i\theta}Z_a^{\theta}$ на угол $-\theta$  является распределением точек $e^{-i\theta}Z_a^{\theta}$  конечной внешней плотности Редхеффера.

\begin{example}
Если при некотором $a\in (0,\pi/2]$  для распределения точек $Z_a^{\theta}$ из  \eqref{Zat} сумма 
\begin{equation*}
\sum_{z\in Z_a^{\theta}}\frac{1}{|z|}<+\infty
\end{equation*}
конечна, то   при выборе $c:=0$ ряд  \eqref{RdZ}, очевидно, сходится, откуда 
распределение точек $Z$ на $\CC$ конечной  внешней плотности Редхеффера вблизи   направления\/ $\theta$. 
Согласно примеру \ref{ex2}, если  для распределения точек  $Z_a^{\theta}$ из \eqref{Zat} при некотором $a\in (0,\pi/2]$ его поворот $e^{-i\theta}Z_a^{\theta}$  на угол $-\theta$  разделённый и  полностью лежит  в какой-нибудь полосе \eqref{strip}  конечной	 ширины $b\in \RR^+$, то $Z$ конечной внешней плотности Редхеффера вблизи  направления\/ $\theta$.
\end{example}

Следующая, более тонкий, чем теорема \ref{th1},  результат значительно обобщает  критерий полноты Мальявена\,--\,Рубела и даже  неулучшаем, если использовать в критерии 
лишь логарифмическую субмеру распределения точек $Z$ на $\CC$.

\begin{thm}\label{th2} Для любых числа $b\in \RR^+$, направления $\theta\in \RR$ и распределения точек $Z$ на $\CC$  конечной внешней плотности Редхеффера вблизи противоположных направлений\/ $\theta$ и $\theta-\pi$ следующие четыре  утверждения равносильны:
\begin{enumerate}[{\rm I.}]
\item\label{KhI0}  Для любого выпуклого компакта $K$ ширины
$\width_K(\theta)\leq b$ экспоненциальная система $\Exp^Z$ полна в  банаховом пространстве $C(K)\bigcap \Hol(\intr K)$ с нормой \eqref{CKn}.  

\item\label{KhII0-} Логарифмическая субмера 
$\mathcal{L}_{e^{i(\pi/2-\theta)}Z}$ для поворота  $e^{i(\pi/2-\theta)}Z$ распределения точек $Z$  на угол $\pi/2-\theta$   удовлетворяет равенству 
\begin{equation}\label{ellZb}
\sup_{1\leq r<R<+\infty}
\biggl(\mathcal{L}_{e^{i({\pi}/{2}-\theta)}Z}(r,R)-\frac{b}{2\pi}\ln\frac{R}{r}\biggr)=+\infty.
\end{equation}

\item\label{KhII0} 
При  значениях $n$ и $N$, пробегающих   соответственно\/ $\NN_0$ и\/ $\NN$,  имеем 
\begin{equation}\label{ellZMg0}
\limsup\limits_{N\to  \infty}\sup\limits_{0\leq n<N}
\biggl({\mathcal L}_{e^{i(\pi/2-\theta)}Z}(e^n,e^N)-\frac{b}{2\pi}(N-n)\biggr)=+\infty.
\end{equation}

\item\label{KhIII0} 
Система $\Exp^Z$ полна в пространстве  $C\bigl(e^{i(\theta-\pi/2)}\overline \strip_{b/2}\bigr)\bigcap \Hol\bigl(e^{i(\theta-\pi/2)}\strip_{b/2}\bigr)$, снабжённом топологией равномерной  сходимости на компактах из замкнутой полосы $e^{i(\theta-\pi/2)}\overline \strip_{b/2}$.
\end{enumerate} 
\end{thm}

\section{Доказательства критериев}
\begin{proof}[Доказательство теоремы\/ {\rm \ref{th1}.}]
Эквивалентность утверждений \ref{KhI+} и \ref{KhI}, верного для любых  систем   целых функций, а не только экспоненциальных,  --- простое следствие определений полноты в функциональных пространствах и топологии равномерной сходимости на компактах.

Докажем эквивалентность утверждений \ref{KhIII} и \ref{KhI}. Неполнота  системы  $\Exp^{Z}$   в пространстве $\Hol\bigl(e^{i(\theta-\pi/2)}\strip_{b/2}\bigr)$ в топологии равномерной сходимости на компактах  по теореме Хана\,--\,Банаха эквивалентна  существованию  ненулевого линейного непрерывного  функционала на $\Hol\bigl(e^{i(\theta-\pi/2)}\strip_{b/2}\bigr)$, 
который обращается в  нуль на каждой экспоненциальной функции из $\Exp^{Z}$ \cite[теорема~1.1.1]{Khsur}. Такой функционал   может быть продолжен как ненулевой линейный непрерывный  функционал  на $C(K)\bigcap \Hol(\intr K)$
 для некоторого выпуклого компакта $K\subset e^{i(\theta-\pi/2)}\strip_{b/2}$, который по-прежнему равен нулю на каждой экспоненциальной функции из $\Exp^{Z}$. Вновь по теореме Хана\,--\,Банаха  \cite[теорема~1.1.1]{Khsur} это означает, что система  $\Exp^{Z}$   не полна в $C(K)\bigcap \Hol(\intr K)$, в то время  как ширина компакта $K\subset e^{i(\theta-\pi/2)}\strip_{b/2}$ в направлении $\theta$ строго меньше  $b$.

Докажем эквивалентность  утверждений \ref{KhII} и \ref{KhIII}. 
Если $Z$ конечной верхней плотности, т.е.  $\overline \dens(Z)<+\infty$, то левая часть   \eqref{lndZ} ---  это один из вариантов определения  \textit{логарифмической блок-плотности распределения точек\/} $e^{i({\pi}/{2}-\theta)}Z$ \cite{Kha89}, \cite[определения (3.2.4)]{Khsur}, а условие \eqref{lndZ} согласно  \cite[теорема 2]{Kha89}, \cite[теорема 3.2.2]{Khsur} означает, что  система $\Exp^{e^{i({\pi}/{2}-\theta)}Z}$  полна в пространстве $\Hol(\strip_{b/2})$. 
Используя обратный поворот на угол $\theta-{\pi}/{2}$  к распределению точек  $e^{i({\pi}/{2}-\theta)}Z$ получаем полноту системы $\Exp^Z$ в пространстве $\Hol\bigl(e^{i(\theta-\pi/2)}\strip_{b/2}\bigr)$. 
 \end{proof}

\begin{proof}[Доказательство теоремы\/ {\rm \ref{th2}.}] Эквивалентность утверждений  
\ref{KhI0} и \ref{KhIII0} обосновывается практически так же, как доказывалась эквивалентность утверждений 
\ref{KhI} и \ref{KhIII} теоремы \ref{th1} выше на основе следствий из теоремы Хана\,--\,Банаха по схеме, изложенной в \cite[п.~1.1.1.]{Khsur}. 

При  доказательстве эквивалентности \ref{KhII0-}$\Leftrightarrow$\ref{KhIII0}
 после соответствующих поворотов $\CC$  и распределения точек $Z$  достаточно рассмотреть случай $\theta =\pi/2$. Выделим в этом случае из $Z$ часть 
$Z_a^{\pi/2}$, определённую в \eqref{Zat},   конечной внешней плотности Редхеффера вблизи $\pi/2$. 
Это означает, что  распределение точек $e^{-i\pi/2}Z_a^{\pi/2}$ конечной внешней плотности Редхеффера, 
из определения которой в соответствии с \eqref{RdZ} легко видеть, что конечна сумма 
\begin{equation*}
\sum_{z\in Z_a^{\pi/2}}\Bigl|\Re\frac{1}{z}\Bigr|<+\infty.
\end{equation*}
Точно так же для части $Z_a^{-\pi/2}$, определённой в \eqref{Zat}, имеем 
\begin{equation*}
\sum_{z\in Z_a^{-\pi/2}}\Bigl|\Re\frac{1}{z}\Bigr|<+\infty.
\end{equation*}
Таким образом, для распределения точек $Z_a^{\pm\pi/2}:=Z_a^{\pi/2}\bigcup Z_a^{-\pi/2}$ сходится сумма 
\begin{equation}\label{Zpmp}
\sum_{z\in Z_a^{\pm\pi/2}}\Bigl|\Re\frac{1}{z}\Bigr|<+\infty.
\end{equation}
По построению часть $Z\setminus Z_a^{\pm\pi/2}$ распределения точек 
отделена от мнимой оси $i\RR$ в том смысле, что она 
лежит вне пары открытых вертикальных углов, содержащих $i\RR\setminus 0$, и при выполнении 
\eqref{ellZb} и \eqref{Zpmp} точная верхняя грань 
\begin{equation}\label{ellZ0}
\sup_{1\leq r<R<+\infty}
\biggl(\mathcal{L}_{Z\setminus Z_a^{\pm\pi/2}}(r,R)-\frac{b}{2\pi}\ln\frac{R}{r}\biggr)
\end{equation}
равна $+\infty$. Отсюда для отделённого от  мнимой оси распределения точек  $Z\setminus Z_a^{\pm\pi/2}$
следует \cite[следствие 4.2]{Kha91}, что экспоненциальная система 
$\Exp^{Z\setminus Z_a^{\pm\pi/2}}$ полна в пространстве 
$C\bigl(\overline \strip_{b/2}\bigr)\bigcap \Hol\bigl(\strip_{b/2}\bigr)$, снабжённом топологией равномерной  сходимости на компактах. Тем более экспоненциальная система   $\Exp^{Z}$ с  б\'ольшим   распределением показателей $Z\supset Z\setminus Z_a^{\pm\pi/2}$ полна в этом пространстве. Таким образом, доказана импликация  \ref{KhII0-}$\Rightarrow$\ref{KhIII0}. 

Теперь из отрицания утверждения \ref{KhII0-} выведем отрицание утверждения \ref{KhIII0}. 
Если не выполнено  \eqref{ellZb} при $\theta=\pi/2$, то   ввиду \eqref{Zpmp} 
конечна точная верхняя грань из \eqref{ellZ0}. Рассмотрим целую функцию 
\begin{equation}\label{gs}
g\colon z\underset{z\in \CC}{\longmapsto} \sin \frac{bz}{2} 
\end{equation}
экспоненциального типа. Легко подсчитать, что найдётся число $C_b\in \RR^+$,  с  которым 
\begin{equation*}
\frac{b}{2\pi}\ln\frac{R}{r}\leq \frac{1}{2\pi}\int_r^R\frac{\ln \bigl|g(iy)g(-iy)\bigr|}{y^2}\dd y+C_b
\quad\text{при всех  $1\leq r<R<+\infty$.}
\end{equation*}
Исходя из последнего согласно предполагаемой конечности  \eqref{ellZ0} для функции $g$ из \eqref{gs} имеем 
\begin{equation*}
\sup_{1\leq r<R<+\infty}
\biggl(\mathcal{L}_{Z\setminus Z_a^{\pm\pi/2}}(r,R)-\frac{1}{2\pi}\int_r^R\frac{\ln \bigl|g(iy)g(-iy)\bigr|}{y^2}\dd y\biggr)<+\infty.
\end{equation*}
При этом условии для целой функции $g$ экспоненциального типа с распределением нулей, лежащим вне пары открытых вертикальных углов, содержащих мнимую ось, а точнее, для \eqref{gs} на  $\RR$, 
найдётся \cite[основная теорема]{Kha91} такая целая функций $f_b\neq 0$ экспоненциального типа с распределением нулей $\Zero_{f_b}$, рассматриваемым с учётом кратности, что $Z\setminus Z_a^{\pm\pi/2}\subset \Zero_{f_b}$ и в то же время $\bigl|f_b(iy)\bigr|\leq \bigl|g(iy)\bigr|$ при всех $y\in \RR$. 
Учитывая явный вид функции $g$ из \eqref{gs}, получаем 
\begin{equation}\label{fbz}
\ln \bigl|f_b(iy)\bigr|\leq \frac{b}{2}|y|+c_b\quad\text{для некоторого $c_b\in \RR^+$ при всех $y\in \RR$,}
\qquad Z\setminus Z_a^{\pm\pi/2}\subset \Zero_{f_b}.
\end{equation}  
Поворот $e^{-i\pi/2}Z_a^{\pm\pi/2}$ распределения точек $Z_a^{\pm\pi/2}$ на угол $-\pi/2$ 
по условию конечной внешней плотности Редхеффера. Для такого распределения точек из сочетания 
теоремы Бёрлинга\,--\,Мальявена о радиусе полноты  с классической теоремой Пэли\,--\,Винера
(см., например,  \cite{Red77}, \cite{Kha94}) следует, что найдётся целая функция $q\neq 0$ 
экспоненциального типа с распределением корней $\Zero_q\supset e^{-i\pi/2}Z_a^{\pm\pi/2}$, для которой  
\begin{equation*}
\bigl|q(x)\bigr|\leq \frac{1}{1+x^2} \quad\text{при всех $x\in \RR$.}
\end{equation*}
Обратный поворот на угол $\pi/2$  даёт  определяемую тождеством 
$h(z)\underset{z\in \CC}{\equiv} q(e^{-i\pi/2}z)$ целую функцию $h\neq 0$ экспоненциального типа, для которой  
\begin{equation*}
\bigl|h(iy)\bigr|\leq \frac{1}{1+y^2} \quad\text{при всех $y\in \RR$,}
\quad   Z_a^{\pm\pi/2}\subset \Zero_h.
\end{equation*}
Отсюда согласно \eqref{fbz}  целая функция $f:=f_bh\neq 0$ экспоненциального типа для некоторой постоянной $C\in \RR^+$ удовлетворяет условиям 
\begin{equation*}
\bigl| f(iy)\bigr|\leq \frac{C}{1+y^2}\exp\frac{b|y|}{2}\quad\text{при всех $y\in \RR$,}
\quad Z\subset \Zero_f. 
\end{equation*}
Это значит \cite[гл.~IV, \S~1, п.~3]{Leo}, что  найдётся прямоугольник $Q\subset \strip_{b/2}$ и ненулевая непрерывная функция $p$ на его границе $\partial Q$, для которых имеет место представление функции $f$ через преобразование Фурье\,--\,Лапласа 
\begin{equation*}
f(w)=\int_{\partial Q} e^{zw}p(z)\dd z.
\end{equation*}
Таким образом, найден линейный непрерывный функционал на $C(\overline \strip_{b/2})\bigcap \Hol(\strip_{b/2})$
\cite[гл.~1 и п.~3.2.2]{Khsur}, который равен нулю на каждой функции из  экспоненциальную систему  $\Exp^Z$. Следовательно, 
система $\Exp^Z$ не полна в  $C(\overline \strip_{b/2})\bigcap \Hol(\strip_{b/2})$ и 
импликация \ref{KhIII0}$\Rightarrow$\ref{KhII0-} тоже доказана. 

Импликация \ref{KhII0-}$\Rightarrow$\ref{KhII0} получается при выборе значений $r=e^n$ и $R=e^N$, 
когда  $n$ и $N$ пробегают   соответственно\/ $\NN_0$ и\/ $\NN$.
Если же \ref{KhII0-} неверно и левая часть \eqref{ellZb} конечна,  то, переходя к случаю $\theta=\pi/2$ и разбирая рассмотренные выше распределения точек $Z_a^{\pm\pi/2}$ и $Z\setminus Z_a^{\pm\pi/2}$, нетрудно убедиться, что каждое из этих распределений точек конечной верхней плотности. Следовательно, их объединение $Z$
тоже конечной верхней плотности и тогда величины ${\mathcal L}_Z(e^n, e^{n+1})$ согласно 
 определениям \eqref{ell}--\eqref{leZ} равномерно ограничены по $n\in \NN_0$.  Отсюда 
при конечности  левой части \eqref{ellZb} нетрудно удостовериться в конечности 
левой части в \eqref{ellZMg0}, что даёт импликацию \ref{KhII0}$\Rightarrow$\ref{KhII0-}. 
\end{proof}

\begin{remark}
Более разнообразные результаты могут быть получены из  электронной публикации \cite{Kha22A}, принятой к печати в существенно расширенной форме в журнал <<Известия РАН. Серия математическая>> \cite{Kha23}. Частично они   доложены на Международных  конференциях по комплексному анализу памяти А.А. Гончара и А.Г. Витушкина  в Москве  \cite{Kha22D} в 2021 г., 
 а также на Международных конференциях <<Понтрягинские чтения -- XXXIII>> в  Воронеже \cite{KhaMur22}
 и <<Комплексный анализ и смежные вопросы>> в Казани 
\cite{Kha22KFU} в 2022 г.

\end{remark}

\section{Условия полноты в иных терминах, выражаемых через ширину} 

\subsection{Широта и диаметр через ширину}
{\it Широтой\/} \cite[гл.~I, \S~4]{Sant}, \cite[п.~3.2 ]{Kha14II}, или {\it толщиной\/} \cite[4.1.1]{Had},
множества $S$  называется наименьшее значение его  ширины  по направлениям 
\begin{equation}\label{breadth}
\breadth (S):=\inf_{\theta \in \RR} \width_S(\theta)\in \overline \RR^+\bigcup -\infty,
\end{equation}
равное   $-\infty$ тогда и только тогда, когда  $S$ пусто, и равное $+\infty$, если и только если  нет полосы конечной ширины,  содержащей $S$.    \textit{Диаметр\/} $\diam (S):=\sup\bigl\{|z-w|\bigm| z\in S, w\in S\bigr\}$ равен наибольшему значению его ширины по направлениям \cite[гл.~I, \S~4]{Sant}, \cite[\S~11, предложение 11.1]{Kha09II} 
\begin{equation}\label{diamw}
\diam (S)=\sup_{\theta \in \RR} \width_S(\theta)\in \overline \RR^+\bigcup -\infty,
\end{equation}
Возможно определение ширины множества $S\subset \CC$ в направлении и через его
\textit{опорную функцию,\/} которую  по традиции \cite{Levin56}, \cite{Levin96}  можно  
определить как  $2\pi$-периодическую функцию
\begin{equation}\label{spf}
{\spf}_S\colon \theta\underset{\theta \in \RR}{\longmapsto} \sup_{s\in S}\Re se^{-i\theta}\in \overline \RR
\end{equation} 
на $\RR$, хотя в теории выпуклости и её применениях более уместно опорной функцией множества $S\subset \CC$ называть положительно однородную  выпуклую функцию \cite{BF}, \cite{Roc}, \cite{Leicht}, \cite{Tikh}, \cite{MI-Tikh}, \cite[п.~1.1]{Kha14II}
\begin{equation}\label{Spf}
{\Spf}_S\colon z\underset{z\in \CC}{\longmapsto}\sup_{s\in S} \Re s\Bar z\underset{z\in \CC}{\overset{\eqref{spf}}{\equiv}}
|z|{\spf}_S(\arg z), 
\end{equation}
где для последнего равенства и всюду далее используем соглашение $0\cdot x=0$ для любого $x\in \overline \RR$.
В терминах опорных прямых \cite[{\bf 33}]{BF}, \cite[4.1.1]{Had}, \cite[гл.~I, \S~4]{Sant}
к  $S\subset \CC$ ширина $\width_S(\theta)$  в направлении $\theta\in \RR$ равна   расстоянию между двумя опорными прямыми к $S$, ортогональными радиус-вектору точки $e^{i\theta}$, что в терминах опорных  функций  \eqref{spf}--\eqref{Spf} означает тождества
\begin{equation}\label{wsp}
\width_S(\theta)\underset{\theta\in \RR}{\equiv}\spf_S(\theta)+\spf_S(\theta+\pi)\overset{\eqref{Spf}}{\underset{\theta\in \RR}{\equiv}}
\frac{1}{r}\bigl(\Spf_S(re^{i\theta})+\Spf_S(-re^{i\theta})\bigr)\quad\text{при всех $r\in \RR^+\setminus 0$}.
\end{equation}
Таким образом, согласно \eqref{wsp} широта  и диаметр  множества $S\subset \CC$ 
 определяются через опорные  функции $\spf_S$ и $\Spf_S$ из \eqref{spf} и \eqref{Spf} по формулам 
\begin{align}\label{breadths}
\breadth (S)&\overset{\eqref{breadth}}{=}\inf_{\theta \in \RR} \bigl(\spf_S(\theta)+\spf_S(\theta+\pi)\bigr)
\overset{\eqref{breadth}}{=}\inf_{0\neq z\in \CC}\frac{\Spf_S(z)+\Spf_S(-z)}{|z|},\\
\diam (S)&\overset{\eqref{diamw}}{=}\sup_{\theta \in \RR} 
\bigl(\spf_S(\theta)+\spf_S(\theta+\pi)\bigr)
\overset{\eqref{diamw}}{=}\sup_{0\neq z\in \CC}\frac{\Spf_S(z)+\Spf_S(-z)}{|z|},
\label{diamsf}
\end{align}

\subsection{Логарифмическая блок-плотность и следствия из теоремы \ref{th1}}
Пусть  $\ell$ ---  функция интервалов $(r,R]\subset \RR^+$ со значениями  из  $\overline\RR$ и $\ell (r,R):=\ell \bigl((r,R]\bigr)$.  Определим окончательно введённые в  \cite{KKh00} после работ \cite{Kha88} и \cite[формула (1.9)]{Kha89}, развивающих \cite[определения 3.4, 3.5]{MR},   четыре \textit{логарифмические блок-плотности\/} \cite[определение 4]{SalKha20}:
\begin{subequations}\label{dens}
\begin{align}
\ln\text{\!-}\overline\dens(\ell )&:=\limsup_{a\to +\infty} \frac{1}{\ln a} \limsup_{r\to+\infty} \ell (r,ar); 
\tag{\ref{dens}$^-$}\label{{dens}barl}\\
\ln\text{\!-\!}\underline{\dens}(\ell )&:=\liminf_{a\to +\infty} \frac{1}{\ln a} \limsup_{r\to+\infty} \ell (r,ar); 
\tag{\ref{dens}$_-$}\label{{dens}l}\\
\ln\text{\!-\!}\dens_{\inf}(\ell )&:=\inf_{a>1} \frac{1}{\ln a} \limsup_{r\to+\infty} \ell (r,ar); 
\tag{\ref{dens}i}\label{{dens}infl}\\
\ln\text{\!-\!}\dens_{\rm b}(\ell )&:=\inf
\left\{b\in \RR^+\colon \sup_{0< r<R<+\infty}\left(\ell (r,R)-b\ln\frac{R}{r}\right) <+\infty\right\}.
\tag{\ref{dens}b}\label{{dens}bl}
\end{align}
\end{subequations}
Положительную функцию интервалов $\ell \geq 0$ называем \textit{логарифмической субмерой интервалов,\/}  если выполнены следующие два условия:
\begin{enumerate}[{\rm [{\bf l}1]}]
\item\label{l1} $\sup\limits_{0< r\in \RR^+} \ell (r,2r)<+\infty$ {\it (логарифмический рост)};
\item\label{l2} $\ell (r_1, r_3)\leq \ell (r_1, r_2)+\ell (r_2, r_3)$ для всех\/ $0< r_1<r_2<r_3<+\infty$ {\it (субаддитивность)}. 
\end{enumerate}
Если в первом неравенстве из\/ {\rm [{\bf l}\ref{l2}]}  знак неравенства $\leq$ можно заменить на знак равенства $=$ для любых значений\/ $0< r_1<r_2<r_3<+\infty$ {\it (аддитивность),} то функцию интервалов $l$ называем  {\it логарифмической мерой интервалов.\/} Очевидно,  класс всех логарифмических субмер интервалов --- выпуклый конус над $\RR^+$, замкнутый относительно операции максимума.

\begin{proposition}[{\rm \cite[теорема 1]{KKh00}, \cite[\S~1]{Kha89}, \cite[предложение 6]{SalKha20}}]\label{clb}
Для логарифмической субмеры интервалов $\ell \geq 0$ все четыре логарифмические блок-плотности из  \eqref{dens} конечны и совпадают, а верхний предел $\limsup\limits_{a\to +\infty}$ в  \eqref{{dens}barl} и нижний предел   $\liminf\limits_{a\to +\infty}$ в \eqref{{dens}l} можно заменить на  обычный предел $\lim\limits_{a\to +\infty}$. Далее для логарифмической субмеры интервалов $\ell \geq 0$  все четыре логарифмические блок-плотности из  \eqref{dens} обозначаем единообразно как $\ln\text{\!-\!}\dens (\ell )$.
\end{proposition} 

\begin{example}
Для  распределения точек $Z$ на $\CC$  \textit{конечной верхней плотности\/}
  примерами логарифмической меры и субмеры интервалов могут служить соответственно 
 определённая в \eqref{ell} правая логарифмическая мера  $\ell_{Z}$  для этого распределения точек $Z$,
 а также логарифмическая субмера  $\mathcal{L}_{Z}$ из \eqref{leZ}. В таком случае 
логарифмическую блок-плотность $\ln\text{\!-\!}\dens (\mathcal{L}_{Z})$ будем обозначать через 
$\ln\text{\!-\!}\dens (Z)$ и называть её {\it логарифмической блок-плотностью распределения точек\/} $Z$.
В частности, левая часть в \eqref{lndZ}  --- это в точности логарифмическая блок-плотность 
 поворота  $e^{i(\pi/2-\theta)}Z$ распределения точек $Z$  на угол $\pi/2-\theta$.
\end{example}

Следующие два  следствия теоремы \ref{th1} иллюстрируют  возможности её применения в случае геометрических характеристик области, отличных от ширины в направлении, хотя и тесно связанных с ней через формулы \eqref{breadth}, \eqref{breadths}, \eqref{diamw}, \eqref{diamsf}.

\begin{corollary}\label{th3} Пусть $Z$ ---  распределение точек на $\CC$ конечной верхней плотности, а также   $0<b\in \RR^+$. 
Тогда следующие четыре  утверждения равносильны: 
\begin{enumerate}[{\rm I.}]

\item\label{espI0}  Для любой выпуклой области $D\subset \CC$ широты $\breadth (D)\leq b$
cистема ${\Exp}^Z$ полна в  $\Hol(D)$.

\item\label{espI} Для любого выпуклого компакта $K$ широты $\breadth (K)< b$
cистема ${\Exp}^Z$ полна в пространстве  $C(K)\bigcap \Hol(\intr K)$.  

\item\label{espII} Выполнено неравенство $\inf\limits_{\theta\in \RR}\ln\text{\!-\!}\dens (e^{i\theta}Z)\geq {b}/{2\pi}$.

\item\label{espIII} Для любого $\theta\in \RR$ система ${\Exp}^Z$ полна в  $\Hol\bigl(e^{i\theta}\strip_{b/2}\bigr)$.
\end{enumerate}
 \end{corollary}
\begin{proof} По определениям широты \eqref{breadth} и \eqref{breadths} эквивалентности теоремы \ref{th3}  сразу следуют из соответствующих эквивалентностей  теоремы \ref{th1}, применённых по всем  $\theta\in \RR$.
\end{proof}
Следующий результат   даёт лишь достаточные условия полноты в терминах диаметра.

\begin{corollary}\label{th4} Пусть $Z$ ---  распределение точек на $\CC$ и   $0<b\in \RR^+$. 

Если $\overline \dens(Z)=+\infty $ или  
$ \sup\limits_{\theta\in \RR}\ln\text{\!-\!}\dens (e^{i\theta}Z)\geq {b}/{2\pi}$,
то экспоненциальная система $\Exp^Z$ полна  в пространстве     $C(K)\bigcap \Hol(\intr K)$ для любого выпуклого компакта $K$ диаметра $\diam(K)< b$ и в $\Hol(D)$ для всякой выпуклой области $D$ диаметра $\diam(D)\leq b$.
 \end{corollary}
\begin{proof} По формуле \eqref{diamw} или \eqref{diamsf} для диаметра  это достаточное условие полноты  системы  $\Exp^Z$ в $C(K)\bigcap \Hol(\intr K)$ для любого выпуклого компакта $K$ диаметра $\diam(K)< b$ 
следует из импликации \ref{KhII}$\Rightarrow$\ref{KhI} теоремы \ref{th1}, что влечёт за собой  полноту  $\Exp^Z$
в $\Hol(D)$ для всех выпуклых областей $D$ диаметра $\diam(D)\leq b$ по определению топологии равномерной сходимости на компактах в $\Hol(D)$. 
\end{proof}
\begin{remark} Одна  наша конструкция \cite[\S~5, теорема 6]{Kha09}
 позволяет для любого сколь угодно большого числа $b>0$  построить разделённое распределение попарно различных точек $Z$ на положительной полуоси $\RR^+$,  для которой 
$\sup\limits_{\theta\in \RR}\ln\text{\!-\!}\dens (e^{i\theta}Z)=0$
и в то же время система $\Exp^{Z}$ полна в любом 
 $C(K)\bigcap \Hol(\intr K)$  при  $\diam (K)<b$. Таким образом, логарифмическая блок-плотность даже разделённого распределения точек на луче    не может полностью характеризовать полноту экспоненциальной системы в терминах диаметра $\diam (K)$. 
\end{remark}

\end{document}